\documentclass{amsart}
\usepackage{graphicx,,amssymb}
\usepackage{amsfonts}
\usepackage{mathrsfs}
\usepackage[all]{xy}
\usepackage{overpic}
\newtheorem{thm}{Theorem}[section]
\newtheorem{lemma}[thm]{Lemma}
\newtheorem{prop}[thm]{Proposition}
\newtheorem{cor}[thm]{Corollary}
\newtheorem{conjecture}[thm]{Conjecture}

\theoremstyle{definition}
\newtheorem{defi}[thm]{Definition}
\newtheorem{example}[thm]{Example}

\theoremstyle{remark}
\newtheorem{remark}[thm]{Remark}
\numberwithin{equation}{section}

\newcommand{\<}{\langle}
\renewcommand{\>}{\rangle}
\newcommand{\st}{\; | \;}             

\newcommand{\PP}{\mathbb{P}}      
\newcommand{\OO}{\mathcal{O}}     
\newcommand{\F}{\mathcal{F}}      
\newcommand{\D}{\mathcal{D}}      
\newcommand{\C}{\mathcal{C}}      
\newcommand{\Ihat}{{\widehat{I}}}   

\newcommand{\Om}{\Omega}
\newcommand{\al}{\alpha}
\newcommand{\be}{\beta}
\newcommand{\Ph}{\Phi}

\newcommand{\CC}{\mathbb{C}}       
\newcommand{\Z}{\mathbb{Z}}       

\DeclareMathOperator{\Ind}{Ind}
\DeclareMathOperator{\Hom}{Hom}
\DeclareMathOperator{\RHom}{RHom}
\DeclareMathOperator{\Ext}{Ext}

\DeclareMathOperator{\id}{id}
\DeclareMathOperator{\Path}{Path}
\DeclareMathOperator{\Rep}{Rep}
\DeclareMathOperator{\rk}{rank}
\begin{document}
\title{Coxeter Elements and Periodic Auslander--Reiten Quiver}
\author{A. Kirillov}
 \address{Department of Mathematics, SUNY at Stony Brook, 
            Stony Brook, NY 11794, USA}
    \email{kirillov@math.sunysb.edu}
    \urladdr{http://www.math.sunysb.edu/\textasciitilde kirillov/} 
\author{J. Thind}
 \address{Department of Mathematics, SUNY at Stony Brook, 
            Stony Brook, NY 11794, USA}
    \email{jthind@math.sunysb.edu}

\maketitle

\section*{Introduction}
Traditionally, to study  a root system $R$ one starts by choosing a
set of simple roots $\Pi\subset R$ (or equivalently, polarization of
the root system into positive and negative parts) which is then used
in all constructions and proofs. To prove that the resulting answers are
independent of the choice of simple roots, one uses the fact that
any two choices of simple roots are conjugate under the action of the
Weyl group. For example, this is the usual way to construct the Dynkin
diagram.

In this paper, we discuss a different approach, starting not with a
set of simple roots but with a choice of a Coxeter element $C$ in the Weyl
group. This is motivated by the geometric construction of root systems of ADE
type from McKay correspondence, outlined in the next section; however,
all our proofs are independent of the geometric construction. In
particular, we show that for a simply-laced root system a choice of
$C$ gives rise to a natural construction of the Dynkin diagram, in
which vertices of the diagram correspond to $C$-orbits in $R$; moreover,
it gives an identification of $R$ with a certain subset 
$\Ihat\subset I\times\Z_{2h}$, where $h$ is the Coxeter number. The set
$\Ihat$ has a natural quiver structure; we call it the {\em periodic
Auslander-Reiten quiver}. Each simple root system $\Pi$  compatible with
$C$ gives rise to an orientation $\Om$ of $I$ and to  a subset of $\Ihat$,
corresponding to roots positive with respect to $\Pi$; such a subset can be
identified with the usual Auslander--Reiten quiver of  the
category $\Rep (I,\Om)$. Using this, we give a description of all
possible simple root systems $\Pi$ compatible with a given Coxeter
element.

We also show that $C$ naturally gives rise to a ``desymmetrization'' of
the bilinear form $( \cdot , \cdot )$ on $R$: it gives us a non-degenerate
integral {\bf non-symmetric} bilinear form $\< \cdot ,\cdot \>$ such that
$(x,y)=\<x,y\>+\<y,x\>$. This form is an analog of the
non-symmetric Euler form in the
theory of quiver representations. 

This gives a combinatorial construction of the root system associated
with the Dynkin diagram $I$: roots are vertices of $\Ihat$, and the
root lattice and the inner product admit an explicit description in
terms of $\Ihat$. 

Finally, we relate our construction with the theory of quiver representations. In particular, we relate the subset of $\Ihat$ corresponding to positive roots to the usual Aulander-Reiten quiver, and we identify the non-symmetric bilinear form $\< \cdot ,\cdot \>$ with the Euler form.

Many of the ideas developed in this paper were first suggested in
Ocneanu's work \cite{ocneanu}; however, to the best of our knowledge, his
constructions were never published.

Throughout the paper, we use the following notation:

$R\subset E$---a reduced, irreducible simply-laced root system in a
Euclidean space $E$ with inner product $(\cdot ,\cdot)$

$W$--- Weyl group of $R$

$C\in W$ --- a Coxeter element

$h$---Coxeter number

{\bf Acknowledgements.}
A.K. would like to thank V. Ostrik, A. Ocneanu, Y. Li, and X. He for
helpful discussions. 

\section{Motivation: geometric construction}\label{s:motivation}

In this section, we describe the motivation for the current work,
which comes from the construction of affine root systems from McKay
correspondence. None of the results of this section are used in the
proofs in the remainder  of the paper; however, this section is crucial
for understanding the ideas behind the constructions.

Recall that the classical McKay correspondence is the correspondence
between finite subgroups $G\subset SU(2)$ and affine Dynkin diagrams
of ADE type. Under this correspondence,  vertices of the  Dynkin diagram
correspond to irreducible representations of $G$; we denote this set
of vertices by $I$. It is also possible to construct the
corresponding root system from $G$, without going through the
construction of the Dynkin diagram. One such construction was
suggested in \cite{kirillov}: if we consider the category $\C$ of
$\bar{G}$--equivariant coherent sheaves on $\mathbf{P}^1$ (where
$\bar{G}=G/\{\pm 1 \}$), and the corresponding 2-periodic derived
category $\D=D^b(\C)/T^2$, where $T$ is the translation functor,
then the set of isomorphism classes of indecomposable objects of $\D$ is
an affine root system, with the inner product given by
\begin{equation}\label{e:inner_prod}
\begin{aligned}
&(x,y)=\<x,y\>+\<y,x\>\\
&\<x,y\>=\dim \RHom(x,y)=\dim \Hom(x,y)-\dim\Ext^1(x,y).
\end{aligned}
\end{equation}
Moreover, the set of isomorphism classes of indecomposable objects can
be very explicitly described. In this description of the root system
there is no standard choice of a set of simple roots, but there is a
fixed Coxeter element given by $\F\mapsto\F\otimes\OO(-2).$

However, in this paper we are interested not in affine  but in
finite root systems. For finite root systems of ADE type, there is a
version of McKay correspondence described in \cite{kirillov-ostrik},
\cite{ostrik} (in which $SU(2)$ must be replaced by a quantum group
at root of unity) but there is yet no reasonable construction of the
analog of the category of equivariant sheaves on (quantum)
$\mathbf{P}^1$. However, by analogy with the affine case we make the
following conjecture.

\begin{conjecture}
    Let $I$ be a simply-laced Dynkin diagram with Coxeter number $h$. Then
    there exists a triangulated category $\D$ and an exact  functor
    $\D\to\D\colon \F\mapsto\F(-2)$ \textup{(}``twist''\textup{)} with the
    following     properties:

    \begin{enumerate}
        \item The category $\D$ is 2-periodic: $T^2=\id$
        \item For any $\F\in \D$, we have a canonical functorial isomorphism
            $\F(2h)=\F$, where $h$ is the Coxeter number of $I$.
        \item Let $K$  be the Grothendieck group of $\D$. Define an inner
            product in $K$ by \eqref{e:inner_prod}. Then  the
            set $\Ind\subset K$  of all isomorphism classes of
            indecomposable objects in $\D$  is a simply-laced root
            system, and the map $C\colon[\F]\mapsto[\F(-2)]$ is a Coxeter
            element for  this root system.
        \item Choose a function $p\colon I\to\Z_2$ \textup{(}``parity
            function''\textup{)} so that $p(i)=p(j)+1$ if $i,j$ are
            connected in $I$. Then there is a natural bijection between
            the set $\Ind$ of indecomposable objects in $\D$ and the set
            $$
            \Ihat=\{(i,n)\in I\times\Z_{2h}\st p(i)+n\equiv 0\mod 2\}
            $$
            Under this bijection, the Coxeter element $C$ defined above is
            identified with the operator $(i,n)\mapsto (i,n-2)$. We will
            denote the indecomposable object corresponding 
            to $(i,n)$ by $X_i(n)$.
        \item $\dim\Hom(X_i(n),X_j(n+1))=n_{ij}$
            is the number of edges between $i,j$ in $I$, and \\
            $\dim\Ext^1(X_i(n),X_j(n+1))=0$. More generally, make
            $\Ihat$ into a quiver by connecting, for every $i,j,n$, the
            vertices  $(i,n)$ and $(j,n+1)$ by $n_{ij}$ oriented edges.
            Then  for any  $(i,n), (j,m)\in \Ihat$ we have
            $$
            \Hom(X_i(n), X_j(m)) =\Path((i,n), (j,m))/ J
            $$
            where $\Path$ is the space of paths in the quiver and $J$ is
            some explicitly described quadratic ideal \textup{(}cf. eq.
             (7.4) in         \cite{kirillov}\textup{)}.
       \item Let $h\colon I\to \Z_{2h}$ be a map  satisfying the following
            condition: if $i,j$
            are connected in $I$, then $h(i)=h(j)\pm 1$. We will call any
            such map a ``height function''; it defines an orientation
            $\Om_h$ on $I$ by    $i\to j$ if $h(j)=h(i)+1$.

            For every such lifting $h$, define
            \begin{align*}
            \Ph_h\colon \D&\to \Rep(I,\Om^{op}_h)\\
                                X&\mapsto \bigoplus_{i\in I}
                                     \Hom(X_i(h_i),X)
    \end{align*}
            where $\Om_h^{op}$ is the orientation opposite to $\Om_h$ 
            and $\Rep (I,\Om^{op}_h)$ is the category of representations   
            of corresponding quiver.
            Let $R\Ph_h$ be the corresponding derived functor $\D\to
            \D^b(I,\Om^{op}_h)/T^2$. Then $R\Ph_h$ is an equivalence of
            triangulated categories.
    \end{enumerate}
\end{conjecture}

It should be noted that such a category $\D$ can be constructed from the
representation theory of quivers: by property (6), for any 
orientation $\Omega$, we have an equivalence $\D\simeq D^b(I,\Om)/T^2$,
which can be taken as the definition of $\D$. In this description, the twist
functor $\F\to\F(2)$ is identified with the Auslander--Reiten functor
$\tau$, and all properties of $\D$ can be derived from known results
about quiver representations. 

However, this description has its drawbacks. Most importantly, it
depends on the choice of orientation of the quiver. A better
construction would be independent of any choices.

For affine Dynkin diagrams, an alternative construction of the
category $\D$  using McKay correspondence  was given in
\cite{kirillov}; as mentioned above, in that case
$\D=\D^b(\C)/T^2$, where $\C$ is the category of $G/\{\pm I\}$-equivariant
coherent sheaves on $\mathbf{P}^1$. However, in the affine case the
statement is slightly more complicated, as in this case indecomposable
sheaves include both free sheaves and torsion ones.

For finite Dynkin diagrams, we do not have an analogous description of
the category $\D$. This is the subject of current
research and will be described in forthcoming papers.

This conjecture would immediately imply a number of purely combinatorial
results; for example, it would  imply that the set  ($C$-orbits in $R$)/$C$
is in natural bijection with the vertices of the Dynkin diagram. The
main goal of the  current paper is to give a purely combinatorial
proof of this and other results, not relying on the conjecture above.

\section{Preliminaries}
In this paper $R\subset E$ is a reduced, irreducible, simply laced root
system (type $A$, $D$ or $E$) with invariant symmetric bilinear form
$(\cdot ,\cdot )$ chosen so that $(\alpha ,\alpha )=2$ for
$\alpha \in R$. We denote by $r$ the rank of $R$, and $W$ its Weyl group. 

We will {\bf not} fix a choice of a simple root system. Instead, recall
that for different simple root systems $\Pi,\Pi'$, there is a unique
element $w\in
W$ such that $w(\Pi)=\Pi'$, which therefore gives a canonical bijection
between simple roots $\al\in \Pi$ and $\al'\in\Pi'$. Therefore, it is
possible to use a single index set $I$ for indexing simple roots in each of the
simple roots systems. More formally, this can be stated as follows. 

\begin{prop}
    There is a canonical indexing set $I$, which depends only on the root
    system $R$, such that for any simple root system $\Pi$ there is a
    bijection 
    \begin{align*}
        I&\to\Pi\\
        i&\mapsto \al_i^\Pi
    \end{align*}
    which is compatible with the action of $W$: if $\Pi'=w(\Pi)$, then
    $w(\al^\Pi_i)=\al_i^{\Pi'}$. 
\end{prop}

For any $i,j\in I$, $i\ne j$, we define 
$$
n_{ij}=-(\al_i^\Pi,\al_j^\Pi)\in \{0,1\}
$$
(this obviously does not depend on the choice of simple root system
$\Pi$); taking vertices indexed by $I$ with $i,j$ connected by $n_{ij}$
unoriented edges gives us the Dynkin diagram of $R$; abusing the notation,
we will also denote this diagram by $I$. 

For a given simple root system $\Pi$, we denote by $R_+^\Pi$ the
corresponding set of positive roots. We also denote 
$s_i^\Pi=s_{\al_i^\Pi}$ the  simple reflections. 

\section{Coxeter element and compatible simple root systems}

From now on, we fix a  Coxeter element $C\in W$ and denote by
$h$ the Coxeter number, i.e. the order of $C$. 

\begin{defi}\label{d:compatible}
    We say that a  simple root system  $\Pi=\{\alpha_1,\ldots,\alpha_r\}$
    is {\em compatible} with $C$ if we can write $C=s_{i_1}\cdots s_{i_r}$.
\end{defi}

By definition, for any Coxeter element there exists at least one
compatible simple root system. However,  not every
simple root system is compatible with given $C$. More precisely, we have
the following result.

\begin{lemma}
    For given simple root system $\Pi$  and  Coxeter element $C$
    we have $l^{\Pi}(C) \geq r$, where $l^{\Pi}(C)$ is  the length of a
    reduced expression for $C$ given in terms of the simple
    reflections $s_i^\Pi$.  Moreover,  $C$ is compatible with $\Pi$ if and
    only if  $l^{\Pi}(C)=r$.
\end{lemma}

\begin{proof}
     
    Let $\omega_i\in E$ be the fundamental weights. Then it is
    immediate from the definition  that 
    $s_i(\omega_j)= \omega_j$ for $j\neq i$ and
    $s_i(\omega_i)= -\omega_i + \sum_{j} n_{ij}\omega_{j}$. If
    $l^{\Pi}(C)< r$  then $C= s_{i_1}\cdots s_{i_l}$, and there exists
    $i\in I$ such that $i\neq i_k$ for any $k$. Hence
    $C(\omega_i)=\omega_i$. However the Coxeter element has no fixed
    vectors in $E$ (see \cite[Lemma 8.1]{kostant1}). Thus, $l^{\Pi}(C) \geq r$. 

    Now suppose that $l^{\Pi}(C)=r$, so $C=s_{i_1}\cdots s_{i_r}$. Then the
    argument above shows that  every $i\in I$ must appear in    
    $\{i_{1},\ldots,i_{r}\}$. Since $|I|=r$, it must appear exactly once,
    so $C$ is  compatible with $\Pi$.
\end{proof}

The next proposition describes the set of Coxeter elements
compatible with a fixed set of simple roots.

\begin{prop}\cite{shi} \label{p:orientations}
    
    \begin{enumerate} 
        \item Let  $C$ be a Coxeter element, $\Pi$ a simple root
        system compatible with $C$.  Choose a reduced expression
        for $C$ and define an
        orientation on $I$        as    follows:  
        $\stackrel {i} {\bullet} \to \stackrel {j}{\bullet}$ if
        $n_{ij}=1$ and 
        $i$ precedes $j$ in a reduced expression for $C$: $C=\dots
    s_i^\Pi\dots s_j^\Pi\dots $. Then this orientation does not depend on
        the    choice of a reduced expression for $C$. 

        \item For fixed $\Pi$, the correspondence 
        $$
          \text{\{Coxeter elements compatible with $\Pi$\} }\to \text{\{
         orientations of $I$\} }
       $$
       defined in Part~1, is a bijection. 
        \end{enumerate}
\end{prop}
\begin{proof}
 For a fixed set of simple roots $\Pi$, if we have two expressions
$C=s_{i_1}\cdots s_{i_r}=s_{i_{1}^{\prime}}\cdots
s_{i_{r}^{\prime}}$ for a Coxeter element, we can obtain one from
another by the operations $s_{i}s_{j}\to s_{j}s_{i}$ for $i,j\in I$
satisfying $n_{ij}=0$. Indeed, it is known (see \cite{bourbaki}) that any
two reduced expressions for an element $w\in W$ can be obtained from each
other by using operations $s_is_j\to s_js_i$ if $n_{ij}=0$
and $s_{j}s_{i}s_{j}\to s_{i}s_{j}s_{i}$ if $n_{ij}=1$; since for a Coxeter
element every simple reflection appears only once, the second operation
does not apply. Thus, the orientation on $I$ does not depend on the choice
of reduced expression. 

Conversely, given an  orientation on $I$, we can define a complete order on
$I$: $I=\{i_1,\dots, i_r\}$ so that all arrows are of the form $\stackrel
{i_k} {\bullet} \to \stackrel {i_l}{\bullet}$ with $k<l$; thus, this
orientation is obtained from the Coxeter element $s^\Pi_{i_1}\dots
s_{i_r}^\Pi$. One easily sees that the order is defined uniquely up to
interchanging $i,j$  with $n_{ij}=0$ and thus the Coxeter element is
independent of this choice of order. 
\end{proof}

\begin{example}
For the root system $R=A_{n}$ and
$\Pi=\{\alpha_1=e_{1}-e_{2},\ldots,\alpha_n=e_{n}-e_{n+1}\}$ the
Coxeter element $C=s_{1}\cdots s_{n}$  corresponds to the
orientation $ \stackrel {1} {\bullet} \to \stackrel {2} {\bullet}
\to \cdots \to \stackrel {n} {\bullet}$.
\end{example}

What we are interested in, however, is the opposite direction: given
a fixed Coxeter element $C$, we would like to describe all 
simple root systems  $\Pi$ which are compatible with $C$.

For fixed compatible pair $(C$, $\Pi)$ the following result shows how to
construct another set of simple roots compatible with $C$, and
describes how the corresponding orientations of $I$ relate.

\begin{prop}\label{p:reflection}
    Let $C$, $\Pi$ be compatible and $i\in I$ be a sink \textup{(}or
    source\textup{)} for the corresponding orientation of $I$ as
defined in Proposition~\ref{p:orientations}.
    \begin{enumerate}
    \item
        $C=s_{i_1}\cdots s_{i_{r-1}}s_{i}$  \textup{(}if $i\in I$ is a
        sink\textup{)} or   $C=s_{i}s_{i_1}\cdots s_{i_{r-1}}$ \textup{(}if
        $i\in I$ is a source\textup{)}.
    \item $\Pi'=s^\Pi_{i}(\Pi)$ is also a set of simple roots compatible
        with $C$. In this case, we will say that $\Pi$ is obtained from
        $\Pi'$ by elementary reflection. Note that we only allow applying
        elementary reflection $s_i$ to $\Pi$ when $i$ is a sink or
        source for    orientation defined by $\Pi$. 
    \item The orientation of $I$ corresponding to $\Pi'$ is
        obtained from the orientation corresponding to $\Pi$ by reversing
        the arrows at $i$. Thus a sink becomes a source, and vice versa.
    \end{enumerate}
\end{prop}

\begin{proof}
We will prove the result for $i\in I$ a sink. The proof for a
source is almost identical.

\begin{enumerate}
\item If $i$ is a sink, all $s_j$ which do not commute with $s_i$ must
    precede $s_i$ in  the reduced expression for $C$. Thus, $s_i$ can be
    moved to the end of the reduced expression. 

\item Denoting temporarily $s_j=s_j^\Pi$, $s'_j=s^{\Pi'}_j$, we
see that $s'_j=s_is_js_i$.    Then, using  Part 1, we can write 
    \begin{align*}
        C &= s_{i_1}\cdots s_{i_{r-1}}s_{i} & \\
         &= (s_{i}s_{i})s_{i_{1}}(s_{i}s_{i})s_{i_{2}}(s_{i}s_{i})\cdots
                (s_{i}s_{i})s_{i_{r-1}}s_{i} &\\
         &=s_{i}(s_{i}s_{i_{1}}s_{i})(s_{i}s_{i_{2}}s_{i})\cdots
                (s_{i}s_{r-1}s_{i}) &\\
         &=s_{i}s'_{i_1}\cdots s'_{i_{r-1}}
          =s'_{i}s'_{i_1}\cdots s'_{i_{r-1}}&
    \end{align*}
    Hence $C$ is compatible with $\Pi'$.

\item From $C=s'_{i}s'_{i_1}\cdots s'_{i_{r-1}}$ and   
Proposition~\ref{p:orientations} we see that $i$ is a source for the
    orientation obtained    from $\Pi'$ and that the orientation is
    obtained by reversing    all the arrows to $i$.
    \end{enumerate}
\end{proof}

\begin{thm}\label{t:reflection2}
    Fix a Coxeter element $C\in W$. Then:
    \begin{enumerate}
    \item The map 
      $$
          \text{\{Simple root systems $\Pi$ compatible with $C$\} }\to
            \text{\{orientations of $I$\} }
     $$
     defined in Proposition~\ref{p:orientations} is surjective. Two
    different simple root systems $\Pi,\Pi'$, both compatible with $C$,
    give the same orientation iff $\Pi'=C^k\Pi$ for some $k\in \Z$. 
    
    \item If $\Pi,\Pi'$ are two  simple root systems, both
    compatible with $C$, then $\Pi'$ can be obtained from $\Pi$ by a
    sequence of elementary reflections $s_i$ as in
    Proposition~\ref{p:reflection}. 
    \end{enumerate}
\end{thm}
\begin{proof}
    The fact that the map is surjective easily follows from
    Proposition~\ref{p:reflection} and the fact that any two orientations
    of a graph without cycles can be obtained one from the other by a
    sequence of operations
    $s_i\colon\text{sink}\leftrightarrow\text{source}$. 
    
    Now, assume that two simple root systems $\Pi$, $\Pi'$ give the same
    orientation. Denoting as before $s_i=s_i^\Pi$, $s'_i=s^{\Pi'}_i$, then we
    see that for some complete order of $I$ we have 
$$
C=s_{i_1}\dots s_{i_r}=s'_{i_1}\dots s'_{i_r}
$$
    (note that the order is the same for $s_i$ and $s'_i$!). Let $w\in W$
    be    such that $w(\Pi)=\Pi'$; then $s'_i=w s_iw^{-1}$ and therefore 
    $C=(ws_{i_1}w^{-1})\dots (ws_{i_r}w^{-1})=wCw^{-1}$,
    so $w$ commutes with $C$. However, it is known (\cite{springer}) that
    the centralizer of the Coxeter element is the cyclic group generated by
    $C$. Thus, $w=C^k$. 

    Finally, to prove the last part, note that it is well known  that 
    any  two orientations can be  obtained one from another by a
    sequence of 
    elementary reflections (reversing all arrows at a sink or a source).
    Thus, if $\Pi,\Pi'$ are   compatible with $C$, then applying a
    sequence of elementary reflections $s_i$ as
    in Proposition~\ref{p:reflection}, we can obtain from $\Pi'$ a simple
    root system $\Pi''$ which gives the same orientation as $\Pi$. By
    Part~1,  it means that $\Pi''=C^k\Pi$. But notice that the simple root
    system $C(\Pi)$ can be obtained from $\Pi$ by a sequence of
    reflections $s_i$: namely, if $C=s^\Pi_{i_1}\dots s^\Pi_{i_r}$, then 
    consider the sequence of simple root systems 
    \begin{align*}
	&\Pi_0=\Pi, \quad \Pi_1=s_{i_1}^\Pi(\Pi_0)\\
    &\Pi_2=s_{i_2}^{\Pi_1}(\Pi_1)=s_{i_1}s_{i_2}s_{i_1}(\Pi_1)
          =s_{i_1}s_{i_2}(\Pi)\\
	&\dots\\
        &\Pi_r=s_{i_r}^{\Pi_{r-1}}(\Pi_{r-1})=s_{i_1}\dots
        s_{i_r}(\Pi)=C(\Pi).
    \end{align*}
    One easily sees that $i_k$ is a source for $\Pi_{k-1}$, so the
    above  sequence of elementary reflections is well-defined. 
\end{proof}
\begin{cor}
    For a given Coxeter element $C$, we have a canonical bijection
         $$
          \text{\{Simple root systems $\Pi$ compatible with $C$\} }/C\to
            \text{\{orientations of $I$\} }
         $$
\end{cor}
In particular, this shows that the number of simple root systems compatible
with $C$ is equal to $h2^{r-1}$, where $h$ is the Coxeter number and $r$
is rank. For example, for root system of type $A_{n-1}$, where the Coxeter
number is $n$ and rank is $n-1$, this gives $n2^{n-2}$ (compared with the
number of all simple root systems, equal to $n!$).

We will give a graphical description of the set of all compatible
simple root systems in terms of ``height functions'' later, in
Theorem~\ref{t:height_functions}.

\section{Representatives of $C$-Orbits}
As before, we fix a Coxeter element $C$ and choose  a simple root
system $\Pi$  compatible with $C$. We define an
order  $\leqslant$ on $I$ by $i\leqslant j$ if there exists an
oriented path $i\to \cdots \to j$, with orientation defined by $\Pi$ as
in Proposition~\ref{p:orientations}. In this case, one easily sees that
$s_i$ must precede $s_j$ in the reduced expression for $C$.

Using this relation we define $\beta_{i}^\Pi\in R$ by
\begin{equation}\label{e:beta}
\beta_{i}^\Pi=\sum _{j\leqslant i} \alpha_{i}^\Pi.
\end{equation}

\begin{prop}\label{p:beta}
    \par\indent
    \begin{enumerate}

    \item The $\beta_{i}^\Pi$ are a basis of the root lattice, and 
        $\alpha_{i} =\beta_{i}^\Pi -\sum_{j\to i} \beta_{j}^\Pi$.

    \item Let $C=s^\Pi_{i_{1}}\cdots s^\Pi_{i_{r}}$ be a reduced 
            expression for
            $C$. Then $\beta^\Pi_{i_{k}} =s^\Pi_{i_{1}}\cdots
            {s}^\Pi_{i_{k-1}} (\alpha^\Pi_{i_{k}})$. 

    \item $\{\beta_{1}^\Pi,\ldots ,\beta_{r}^\Pi\}
        =\{\alpha \in R_{+}^\Pi \st C^{-1}\alpha \in R_{-}^\Pi\}$
        where
        $R_{\pm}^\Pi$ are the sets of positive and negative
        roots defined by the simple root system $\Pi$.

    \item $\beta_{i}^\Pi$ are representatives of the $C$-orbits in $R$.
\end{enumerate}
\end{prop}

\begin{remark} In the theory of quiver representations, the simple representations $X(i)$ correspond to the simple roots $\alpha_{i}^{\Pi}$ in $R$, while the projective representations $P(i)$ correspond to the $C$-orbit representatives $\beta_{i}^{\Pi}$. This correspondence will be discussed in detail in Section~\ref{s:identification}.
\end{remark}

\begin{proof}
     The first two  parts are easily obtained by explicit computation (see, e.g.
\cite[Theorem 8.1]{kostant1} ).
        The  other two parts are known; a proof can be found in
\cite[Chapter VI, \S 1, Proposition 33]{bourbaki}.
\end{proof}
\begin{example}
For the root system $R=A_{n}$ with simple roots $\Pi =\{ \alpha_{1}
=e_{1}-e_{2}, \ldots ,\alpha_{n} =e_{n}-e_{n+1} \}$ and Coxeter
element $C=s_{1}s_{2}\cdots s_{n}$ we get $\beta_{i}^\Pi =\sum_{j
\leqslant i} \alpha^\Pi_{j} =e_{1}-e_{i}$.
\end{example}

\begin{example}\label{ex:bipartite}
Let $I=I_{0} \sqcup I_{1}$ be a bipartite splitting and 
$C= (\Pi_{i\in I_{0}} s_{i}) (\Pi_{i\in I_{1}} s_{i})$. Note that in the
corresponding
orientation of $I$ all the arrows go from $I_{0}$ to $I_{1}$. So $I_{0}$
are sources and $I_{1}$ are sinks. Then $\beta_{i}^\Pi =\alpha^\Pi_{i}$
for $i\in I_{0}$ and $\beta_{i}^\Pi =\alpha^\Pi_{i} + \sum_{j}
n_{ij}\alpha^\Pi_{j} =-C(\alpha^\Pi_{i})$ for $i\in I_{1}$. In this case
the $\beta_{i}^\Pi$ we obtain are the same $C$-orbit representatives as
in \cite{kostant} except
that our $\beta_{i}^\Pi$ for $i\in I_{1}$ are shifted by $C$.
\end{example}

\begin{example}\label{ex:beta-D}
As a special case of the previous example, consider the Dynkin diagram
of
type $D_{2n+1}$, with simple root system 
$$
\Pi = \{ \alpha_{1}=e_{1} -e_{2},\dots, \al_{2n-1}=e_{2n-1}-e_{2n},
\al_{2n}=e_{2n}-e_{2n+1}, e_{2n+1}=e_{2n}+e_{2n+1}\}
$$
 and 
 $$
 I_0=\{2,4, \dots, 2n, 2n+1\},\qquad 
 I_1=\{1,3,\dots, 2n-1\}
 $$

The corresponding  Coxeter element $C= (\Pi_{i\in I_{0}} s_{i})
(\Pi_{i\in
I_{1}} s_{i})$, and the representatives of $C$-orbits are 
\begin{align*}
&\beta_{i}^\Pi =\alpha_{i} \text{ for } i\in I_{0}=\{2,4, \dots, 2n,
2n+1\}\\
&\beta_{1}^\Pi =e_{1} -e_{3}, \beta_{3}^\Pi =e_{2} -e_{5},\dots 
e_{2i+1}=e_{2i}-e_{2i+3}, \dots, \beta_{2n-3}=e_{2n-2}-e_{2n-1}\\ 
&\beta_{2n-1}^\Pi =e_{2n-2}+e_{2n}
\end{align*}

Figure~\ref{f:D5} shows the corresponding orientation, roots $\al_i$
and  $\be_i$ for $D_5$. 

\begin{figure}[ht!]
\vspace{5mm}
\begin{overpic}[scale=1.5]%
    {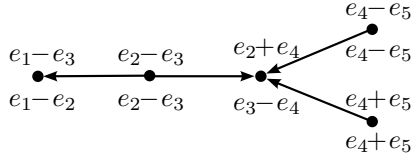}
    \put(-6,8){$e_1\!\!-\!e_2$}
    \put(-6,21){$e_1\!\!-\!e_3$}
    \put(25,8){$e_2\!\!-\!e_3$}
    \put(25,21){$e_2\!\!-\!e_3$}
    \put(58,7){$e_3\!\!-\!e_4$}
    \put(58,23){$e_2\!\!+\!e_4$}
    \put(90,22){$e_4\!\!-\!e_5$}
    \put(90,34){$e_4\!\!-\!e_5$}
    \put(90,-4){$e_4\!\!+\!e_5$}
    \put(90,8){$e_4\!\!+\!e_5$}
\end{overpic}
\caption{The Dynkin Diagram $D_{5}$. $\beta_{i}^\Pi$ is above the
node $i\in I$, and the simple root $\alpha_{i}$ is below the
node $i\in I$.}\label{f:D5}
\end{figure}
\end{example}

We now discuss how the set of $\be_i$ change when we change the simple
root system $\Pi$ (keeping $C$ fixed). By Theorem~\ref{t:reflection2},
it suffices to describe how $\be_i$ change under elementary
reflections. 

\begin{prop}\label{p:elem_reflection_beta}
    Let $\Pi$ be a simple root system compatible with $C$, and let 
    $i\in I$
    be  a sink for the corresponding orientation. Then
    $$
        \beta_j^{s_i\Pi} =\begin{cases} 
                        C^{-1} \beta_{i}^\Pi &\text{ for }j=i \\
                        \beta_{j}^\Pi        &\text{ for } j \neq i
                 \end{cases}
    $$ 
    
    Similarly, if $i$ is a source, then 
    $$
        \beta_j^{s_i\Pi} =\begin{cases} 
                        C \beta_{i}^\Pi &\text{ for }j=i \\
                        \beta_{j}^\Pi   &\text{ for } j \neq i
                 \end{cases}
    $$ 
\end{prop}

\begin{proof} 
    Denote for brevity $s_j=s_j^\Pi$, $s'_j=s_j^{s_i\Pi}$.
    For $i$ a sink, we can write $C=s_{i_{1}}\cdots
    s_{i_{r-1}}s_{i}$ and $\beta_{i} =s_{i_{1}}\cdots s_{i_{r-1}}
    \alpha_{i}$. Hence
    \begin{align*}
    C^{-1}\beta_{i} &= (s_{i}s_{i_{r-1}}\cdots
s_{i_{1}})(s_{i_{1}}\cdots
    s_{i_{r-1}}) \alpha_{i} &\\
    &= s_{i}\alpha_{i} &\\
    &=  \alpha_{i}^{s_{i}\Pi} &\\
    &= \beta_{i}^{s_{i}\Pi} &
    \end{align*}
    since $i$ is a source for $s_{i}\Pi$.\\
    Now take $i_{j} \neq i$. Then
    \begin{align*}
    \beta_{i_{j}}^\Pi &= s_{i_{1}}\cdots s_{i_{j-1}}\alpha_{i_{j}} &\\
    &= (s_{i}s_{i})s_{i_{1}}(s_{i}s_{i})\cdots
    (s_{i}s_{i})s_{i_{j-1}}(s_{i}s_{i})\alpha_{i_{j}} &\\
    &= s_{i}s_{i_{1}}^{\prime}\cdots s_{i_{j-1}}^{\prime}
    (s_{i}\alpha_{i_{j}}) &\\
    &= \beta_{i_{j}}^{s_{i}\Pi} &
    \end{align*}
    Similarly, for $i$ a source we can write $C=s_{i}s_{i_{1}}\cdots
    s_{i_{r-1}}$ and $\beta_{i} =\alpha_{i}$. Hence
    \begin{align*}
    C\alpha_{i} &= s_{i}s_{i_{1}}\cdots s_{i_{r-1}} \alpha_{i} &\\
    &= s_{i_{1}}^{\prime}\cdots s_{i_{r-1}}^{\prime} s_{i}\alpha_{i} &\\
    &= \beta_{i}^{s_{i}\Pi} &
    \end{align*}
    since $i$ is a sink for $s_{i}\Pi$.\\
    For $i_{j}\neq i$ a similar calculation to the case $i$ a sink gives
    us that $\beta_{i_{j}}^\Pi= \beta_{i_{j}}^{s_{i}\Pi}$. 
\end{proof}

\begin{thm}\label{t:C-orbits}
    Let $R/C$ be the set of $C$-orbits in $R$. Then we have a bijection
     \begin{align*}
    I&\to R/C\\
    i&\mapsto C\text{-orbit of }\be_i^\Pi
    \end{align*}
    which does not depend on the choice of a simple root system $\Pi$
    compatible with $C$. 
\end{thm}
\begin{proof}
    The fact that it is a bijection follows from
    Proposition~\ref{p:beta}. To show independence of the choice of
$\Pi$,
        notice that by Proposition~\ref{p:elem_reflection_beta}, if
    $\Pi,\Pi'$ are obtained one from another by an elementary
reflection,
    then the $C$-orbit of $\be_i^\Pi$ and $\be_i^{\Pi'}$ coincide. On the
    other hand, by Theorem~\ref{t:reflection2}, any two simple root
    systems compatible with $C$ can    be obtained one from another by
    elementary reflections. 
\end{proof}
 
For future use, we also give here another proposition, describing the
action of $C$ on $\be_i$. Its geometric meaning will become clear in
Section~\ref{s:root_lattice}.

\begin{prop}\label{p:Cbeta}
 $ C\beta_{i}^\Pi = -\beta_{i}^\Pi+\sum_{j \to i} C\beta_{j}^\Pi +
\sum_{j \leftarrow i} \beta_{j}^\Pi$.
\end{prop}

\begin{proof} 
    Let $i\in I$ be a sink. Then there are no $j\leftarrow i$ and from
    the proof of Proposition~\ref{p:elem_reflection_beta}, 
    $-C^{-1} \beta_{i}^\Pi =\alpha_{i} 
    =\beta_{i}^\Pi -\sum_{j\to i} \beta_{j}^\Pi$. 
    By applying $C$ to this equation and
    rearranging we get the statement of the proposition. 
 
    Now if $i\in I$ is not a sink apply a sequence of reflections
    $s_{j}$  with $j\geq i$ to make it one. This process replaces $\Pi$
    by another  compatible simple root system $\Pi'$. By
    Proposition~\ref{p:elem_reflection_beta},  
    \begin{enumerate}
    \item $\beta_{i}^\Pi =\beta_{i}^{\Pi'}$  
    \item if $j\to i$ for $\Pi$ then $\beta_{j}^\Pi
        =\beta_{j}^{\Pi'}$ 
    \item if $j\leftarrow i$ for $\Pi$ then
            $\beta_{j}^{\Pi'}=C^{-1}\beta_{j}^\Pi.$ 
    \end{enumerate}
    Then by the above calculation we have that 
    \begin{align*}
    \beta_{i}^{\Pi'} + C\beta_{i}^{\Pi'} 
    &= \sum_{j \to i \ \text{in} \ \Pi'} 
        C\beta_{j}^{\Pi'} & \\
    &= \sum_{j \to i \ \text{in} \ \Pi} C\beta_{j}^{\Pi'}  
      +\sum_{j \leftarrow i \ \text{in} \ \Pi} C\beta_{j}^{\Pi'} & \\
    &=\sum_{j \to i \ \text{in} \ \Pi} C\beta_{j}^\Pi 
      + \sum_{j \leftarrow  i \ \text{in} \ \Pi} C(C^{-1}\beta_{j}^\Pi )
      & \\ 
    &=\sum_{j \to i \ \text{in} \ \Pi} C\beta_{j}^\Pi 
      + \sum_{j \leftarrow i \ \text{in} \ \Pi} \beta_{j}^\Pi.& 
    \end{align*}
\end{proof}

\section{Periodic Auslander--Reiten quiver}\label{s:Ihat}
We choose a ``parity function'' $p\colon I\to
\Z_2$ such that $p(i)=p(j)+1$ if $i$ and $j$ are connected by an edge in
$I$. Define
\begin{equation}\label{e:Ihat}
    \Ihat=\{(i,n)\in
    I\times\Z_{2h}\st p(i)+n\equiv 0\mod 2\}.
\end{equation}

We make $\Ihat$ into a quiver by joining $(i,n)$ to $(j,n+1)$ by
$n_{ij}$ oriented edges, where $n_{ij}$ is the number of edges between
$i$ and $j$ in $I$. We will call the resulting quiver the {\em periodic
Auslander--Reiten quiver} and denote it also by $\Ihat$.  
\begin{remark}
    Note that the definition does not require a choice of orientation of
    $I$. However, if we do choose an orientation of $I$, then the quiver
    $\Ihat$ can be identified with  the Auslander--Reiten quiver of the
    2-periodic derived category of representations of quiver $I$. This
    will be discussed in detail in Section~\ref{s:identification} and Section~\ref{s:eulerform}.
\end{remark}

\begin{example}
    For the root system $D_{5}$ the quiver $\Ihat$ is shown in
    Figure~\ref{f:Ihat-D5}. 
    \begin{figure}[ht]
        \includegraphics[height=5.00in]{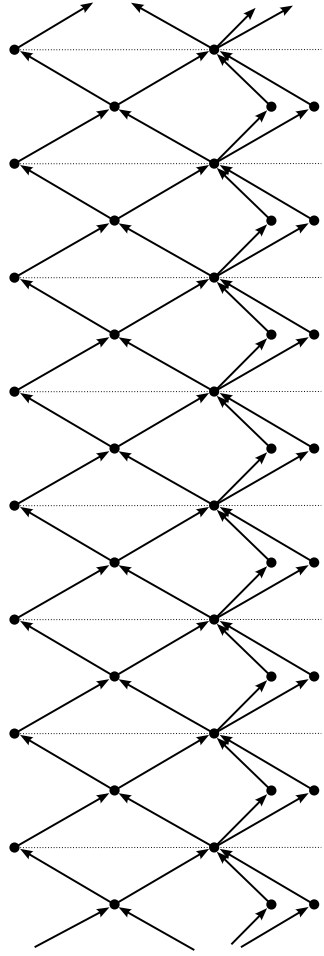}
        \caption{The quiver $\Ihat$ for type $D_{5}$. For comparison, 
        we also show the Dynkin diagram of $D_5$ next to
        it. Note that $\Ihat$ is periodic: the arrows leaving the 
        top level are the same as the arrows going in to the bottom
        level.}\label{f:Ihat-D5}
    \end{figure}
\end{example}

We will define the ``twist'' map  $\tau\colon \Ihat \to
\Ihat$ by
\begin{equation}\label{e:tau}
\tau(i,n)=(i,n+2).
\end{equation}

Our goal is to show that the set of vertices of $\Ihat$ can be (almost)
canonically identified with the root system $R$. To do so, we fix one of
the vertices $i_0\in I$ and choose an identification 
\begin{equation}\label{e:identification}
    \Phi_0\colon (C\text{-- orbit of }\be_{i_0})\to 
        \{(i_0,n)\st n\in \Z_{2h}, n+p(i_0)\equiv 0\mod 2\}
        \subset \Ihat
\end{equation}
which identifies the Coxeter element with the twist:
$\Phi_0(C\be)=\tau\Phi_0(\be)$.  

\begin{remark}
    The seeming arbitrariness in the choice of $\Phi_0$ could have been
    avoided if in the definition of $\Ihat$, instead of $\Z_{2h}$, we had
    used a suitable $\Z_{2h}$---torsor. 
\end{remark}

\begin{figure}[ht]
    \begin{overpic}
    {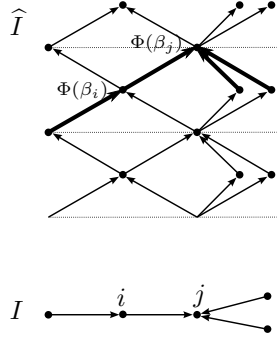}
    \put(4,72){$\scriptstyle{\Phi(\be_i)}$}
    \put(26,86){$\scriptstyle{\Phi(\be_j)}$}
    \put(22,9){$i$}
    \put(45,10){$j$}
    \put(-10,5){$I$}
    \put(-10,90){$\Ihat$}
    \end{overpic}
    \caption{The map $\Phi$. Vertices $\Phi(\be_i)$ and
            edges between them are shown bold.}\label{f:lift}
\end{figure}

\begin{thm}\label{t:Phi1}
    Let $C$ be a fixed Coxeter element and   $\Pi$ a simple root system
    compatible with $C$. Then there    exists a unique bijection
    $\Phi^\Pi\colon R\to\Ihat$ with the following properties
    \begin{enumerate}
        \item It identifies the Coxeter element with the twist: 
        $\Phi^\Pi(C\be)=\tau\Phi^\Pi(\be)$.  

        \item It agrees with the identification $R/C\to I$ given in
            Theorem~\ref{t:C-orbits}:
        $\Phi^\Pi(\be^\Pi_i)=(i,h(i))$ for some $h\colon I\to \Z_{2h}$.
 
        \item If, in the orientation defined by $\Pi$, we have $i\to j$,
            and $\Phi^\Pi(\be_i^\Pi)=(i,n)$, then 
            $\Phi^\Pi(\be_j^\Pi)=(j,n+1)$ 
            \textup{(}see Figure~\ref{f:lift}\textup{)}.
       \item On the $C$--orbit of $\be_{i_0}$, $\Phi^\Pi$ coincides with
        $\Phi_0$.
    \end{enumerate}
\end{thm}
\begin{proof}
    Since $\be_i^\Pi$ are representatives of $C$-orbits in $R$
    (Proposition~\ref{p:beta}), and each
    $C$-orbit has period $h$, it suffices to define
    $\Phi^\Pi(\be_i^\Pi)$.  On the other hand, condition (4) uniquely
    defines $\Phi^\Pi(\be^\Pi_{i_0})$, and it is easy to see that
    given $\Phi^\Pi(\be^\Pi_{i_0})$, conditions (2) and (3) uniquely
    determine $\Phi^\Pi(\be^\Pi_{i})$ for all $i\in I$ (since $I$ is
    connected and simply-connected).  
\end{proof}

\begin{thm}\label{t:Phi2}
    The bijection $\Phi\colon R\to \Ihat$, defined in Theorem~\ref{t:Phi1},
    does  not depend on the choice of $\Pi$.
\end{thm}
\begin{proof}
    By Theorem~\ref{t:reflection2}, it suffices to check that if $\Pi'$
    is obtained from $\Pi$ by  elementary reflection, then
    $\Phi^\Pi=\Phi^{\Pi'}$. To do so, it
    suffices to check that $\Phi^\Pi$ satisfies all the defining
    properties of  $\Phi^{\Pi'}$. The only property which is not obvious
    is (3): if  in the orientation defined by $\Pi'$, we have
    $i\to j$,  and $\Phi^\Pi(\be'_i)=(i,n)$, then 
    $\Phi^\Pi(\be'_j)=(j,n+1)$, where for brevity we denoted 
    $\be_i=\be_i^\Pi$ and $\be'_i=\be^{\Pi'}_i$. 

    Let $\Pi'=s_k\Pi$. If $i,j$ are both distinct from $k$,
    then by Proposition~\ref{p:elem_reflection_beta},    $\be_i=\be'_i$,
    $\be_j=\be'_j$, so  property (3) for $\Phi^{\Pi'}$ coincides with
    the  one for  $\Phi^\Pi$. So we only need to consider the case
    when $i=k$  or $j=k$; in these cases, $s_k$ reverses orientation
    of edge between  $i$ and $j$. 

    If $i=k$,  then for $\Pi$ we have $j\to i$, so $i$ is a sink for
    $\Pi$. Then by Proposition~\ref{p:elem_reflection_beta}
    $\be'_j=\be_j$,  $\be'_i=C^{-1}\be_i$.
     By definition, if $\Phi^\Pi(\be_i)=(i,n)$, then
    $\Phi^\Pi(\be_j)=(j,n-1)$, and 
    $\Phi^\Pi(\be'_i)=\tau^{-1}(i,n)=(i, n-2)$ (see
     Figure~\ref{f:reflections2}). Thus,
    condition  (3) for pair $\be'_i$, $\be'_j$ is satisfied.

\begin{figure}[ht]
    \begin{overpic}
    {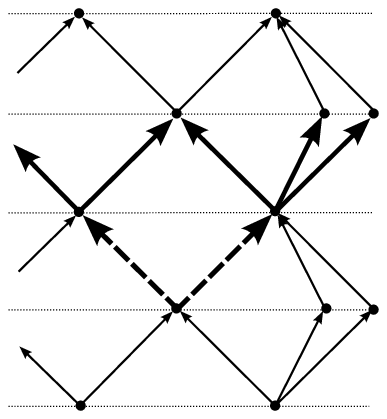}
    \put (37,62){$\Phi(\be_i)$}
    \put (35,20){$\Phi(\be'_i)$}
    \put (81,40){$\Phi(\be_j)$}
    \end{overpic}
    \caption{}\label{f:reflections2}
\end{figure}

    Case $j=k$ is done similarly. 
\end{proof}

\begin{example}\label{ex:phi-A4}
    For Dynkin digram of type $A_4$, with 
    $\Pi=\{e_1-e_2,\dots, e_4-e_5\}$
    and $C=s_{1}s_{2}s_{3}s_{4}$, the map $\Phi\colon
    R\to\Ihat$ is shown in Figure~\ref{f:phi-A4}. 
    For $A_n$ the figure is similar. 
\begin{figure}
\includegraphics{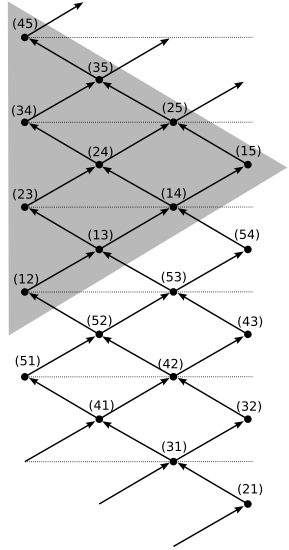}
\caption{Map $\Phi$ for root system of type $A_4$. The figure shows, for
each vertex in $\Ihat$, the corresponding root $\al\in R$; notation
$(ij)$ stands for $e_i-e_j$. The set of positive roots (with respect to
usual polarization of $R$) is shaded. Recall that the quiver $\Ihat$ is periodic.} \label{f:phi-A4}
\end{figure}
\end{example}

As an immediate application of this construction, we will show how this
allows one to give an explicit description of the set of all simple root
systems compatible with $C$.

\begin{defi}\label{d:height}
    A function $h\colon I\to \Z_{2h}$ satisfying $h(j)=h(i)\pm 1$ if
    $i,j$ are connected by an edge in $I$ and satisfying 
    $h(i)\equiv p(i)\mod     2$,  will be called a ``height''
    function. (Here  $p$ is the parity function defined in the beginning
    of this section.)
\end{defi}

\begin{defi}\label{d:slice}
Following \cite{gabriel}, a connected full subquiver of $\Ihat$ which contains a unique representative of $\{ (i,n) \}_{n\in \Z_{2h}}$ for each $i\in I$ will be called a slice.
\end{defi}

Any height function $h$ gives us a slice $I_{h} = \{(i,h(i))\st i\in I\}\subset
\Ihat$; it also defines an  orientation $\Om_h$ on $I$:
$i\to j$ if $i,j$ are connected by an edge and $h(j)=h(i)+1$. It is easy
to see that two height functions give the same orientation iff they
differ by an additive constant, or, equivalently, corresponding slices 
are obtained one form another by applying a power of $\tau$. 

\begin{thm}\label{t:height_functions}
    Let $\Pi$ be a simple root system compatible with $C$. Then the
    function $h^\Pi\colon I\to \Z_{2h}$ defined by $\Phi(\be_i^\Pi)=(i,
    h^\Pi(i))$ is a height function. Conversely, every height function
    can be  obtained in this way from a unique simple root system
    compatible  with $C$.
\end{thm}
\begin{proof}
    The fact that $h^\Pi$ is a height function is obvious from the
    definition of $\Phi$. To check that  any height function can be
    obtained from some $\Pi$, note that by Theorem~\ref{t:reflection2},
    for given height function $h$ we can can find a simple root system
    $\Pi$ which would give the same orientation of $I$ as $h$.
    Therefore, we would have $h^\Pi=h+a$ for some constant $a\in
    \Z_{2h}$, which must  necessarily be even. Take
    $\Pi'=C^{a/2}(\Pi)$; then $h=h^{\Pi'}$. 
\end{proof}
        
\begin{cor}
    For a given Coxeter element $C$, we have a canonical bijection
         $$
          \text{\{Simple root systems $\Pi$ compatible with $C$\} } \to
            \text{\{Height functions on $I$\} }
         $$
\end{cor}

The elementary reflections $s_i$ can also be easily described in terms
of height functions. 

\begin{figure}[ht]
\includegraphics{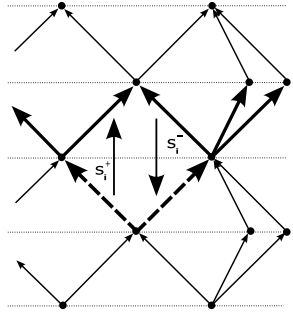}
\caption{Action of simple reflections on height
functions}\label{f:refl_height}
\end{figure}

\begin{prop}
    If $\Pi$ is a simple root system compatible with $C$, $i$---a sink
    for the orientation defined by $\Pi$, and $h^\Pi$, $h^{s_i\Pi}$ ---
    the corresponding height functions as
    defined in Theorem~\ref{t:height_functions}, then
   $$
    h^{s_i\Pi}(j)=\begin{cases} h(j)-2 &j=i \\
                                 h(j) & j \neq i
                 \end{cases}
    $$
    Simiarly, if $i$ is a source, then  
   $$
    h^{s_i\Pi}(j)=\begin{cases} h(j)+2 &j=i \\
                                 h(j) & j \neq i
                 \end{cases}
    $$
	\textup{(}see Figure~\ref{f:refl_height}\textup{)}. 
\end{prop}
The proof is immediate from Proposition~\ref{p:elem_reflection_beta}. 

Abusing the language we will write $s_{i}h$ instead of $h^{s_{i}\Pi}$.

\section{Root lattice}\label{s:root_lattice}

In the previous section, we constructed a canonical bijection
$\Phi\colon R\to \Ihat$, which is independent of the choice of a simple
root system $\Pi$. Among other things, it allowed us to identify the set
of all compatible $\Pi$ with the set of ``height functions'', by using
$\Phi(\be_i^\Pi)$. In this section, we further study this isomorphism.
In particular, we describe the root lattice in terms of $\Ihat$. 
\begin{thm}\label{t:root_lattice}
    Let the lattice $Q$ be defined by 
    \begin{equation}\label{e:Q}
    Q=\Z^{\Ihat}/J,
    \end{equation}
    where $J$ is the ideal generated by the following relations, for
    each $(i,n)\in \Ihat$, 
    \begin{equation}\label{e:fund_relations}
    (i,n)-\sum_{j}(j,n+1)+(i,n+2)=0
    \end{equation}
    \textup{(}the sum is over all vertices $j\in I$ connected to
    $i$\textup{)}. Then 
    \begin{enumerate}
    \item 
        The identification $\Ph\colon R\to \Ihat$ defined in
        Theorem~\ref{t:Phi2}  descends to an isomorphism of lattices
        $Q(R)\to Q$, where $Q(R)$ is the root lattice of $R$.
    \item For any height function $h$, elements $(i,h(i))\in \Ihat$ 
	form a basis of $Q$. 
    \end{enumerate}
\end{thm}
\begin{proof}
      Let $h$ be a height function. Then classes $(i,h(i)),
     i\in I$,  generate the lattice $Q$. Indeed, let $Q_h$ be the
     subgroup generated by $(i,h(i))$. It follows
    from relations  \eqref{e:fund_relations} that if $i$
    is a source for $h$, then $(i,h(i)+2)$ is in $Q_h$;
    thus, $Q_h=Q_{s_i h}$. Since any height function can
    be obtained from $h$ by succesive application of $s_i$,
    we see that $Q_h$ contains all $(i,n)\in \Ihat$, so
    $Q_h=Q$. In particular, this implies that $\rk(Q)\leq  r$. 
  
    Next, we show that $\Ph^{-1}$ descends to a   well-defined map
    $Q\to Q(R)$. To do so, we need to show that 
     $$
      \Ph^{-1}(i,n)-\sum_{j}\Ph^{-1}(j,n+1)+\Ph^{-1}(i,n+2)=0
    $$
    in $Q(R)$. Choosing a simple root system $\Pi$ such that
    $\Ph^{-1}(i,n)=\be_i^\Pi$, we see that this relation is 
    equivalent to the relation 
    $$\beta_{i}^\Pi -\Bigl( \sum_{j \to i} C\beta_{j}^\Pi 
          +\sum_{j \leftarrow i} \beta_{j}^\Pi\Bigr)
    + C\beta_{i}^\Pi =0
    $$
    proved in Proposition~\ref{p:Cbeta}.

    Since $\Phi^{-1}\colon \Ihat\to R$ is a bijection, the map
    $\Phi^{-1}\colon Q\to Q(R)$ is surjective. Since we have already
    shown that $\rk(Q)\leq r$, this implies that $\rk(Q)=r$, so that
    $\Phi^{-1}$ is an isomorphism, and that for fixed height function
    $h$, the classes  $(i,h(i))$ form a basis of $Q$. 
\end{proof}
\begin{remark}
    Relations \eqref{e:fund_relations} are motivated by  almost
    split exact    sequences in the theory of quiver representations, or
    by the short    exact sequence of coherent sheaves on
    $\PP^1=\PP(V)$:    $0\to \F\to \F(1)\otimes V\to \F(2)\to 0$. 
\end{remark}

\section{Euler form}\label{s:euler}

Recall the conjectural picture laid out in
Section~\ref{s:motivation}. In particular, recall the Euler form
$\< X,Y \> =\dim \RHom(X,Y)=\dim \Hom(X,Y)-\dim\Ext^{1}(X,Y)$
defined on $\D$. In this section we define the analog of such a form
on both $R$ and $\Ihat$ and show how it can be used to give  a
definition of the inner product in terms of $\Ihat$.

Define a bilinear form $\< \cdot ,\cdot\>^\Pi$ on the
root lattice $Q(R)$ by
\begin{equation}\label{e:form}
    \< \beta_{i}^\Pi,\alpha_{j}^\Pi \>^\Pi=\delta_{ij}.
 \end{equation}

 By Proposition~\ref{p:beta}, this  completely determines
$\< \cdot ,\cdot \>^\Pi$. Note that this form
is  non-degenerate but not symmetric.

\begin{thm}\label{t:form1}
    \par\indent
    \begin{enumerate}
    \item $\< \beta_{i}^\Pi ,\beta_{j}^\Pi \>^\Pi =
              \begin{cases} 1 & \text{i} \leqslant j \\ 0 &
              \text{otherwise} \end{cases}=$
                  the number of paths $(i\to \cdots \to j)$.

    \item $\< \cdot ,\cdot \>^\Pi$ is integer
        valued on $R$ and satisfies
    
        \begin{align}
            &\< x,y\>^\Pi +
               \< y,x\>^\Pi=(x,y)\label{e:symmetrization}\\
            &\< x,y\>^\Pi =-\< y,C^{-1}
x\>^\Pi\label{e:serre_duality} .
        \end{align}
      where $(\cdot,\cdot)$ is the  $W$-inner product in
    $E$ normalized so that $(\al,\al)=2$ for $\al\in R$. 
    \end{enumerate}
\end{thm}

Note that the equation  $\< x,y\>  =-\< y,C^{-1} x\>$
in $R$ corresponds to the statement of Serre Duality $\Hom(X,Y) =
\Ext^1(Y,X(-2))^{*}$ in the (conjectural) category $\D$ from
Section~\ref{s:motivation}, or equivalently, the identity
$\Ext^1(X,Y)=D \Hom(X,\tau Y)$ in Auslander--Reiten theory. 

\begin{proof}
    Since throughout the proof the simple root system $\Pi$ will be
    fixed, we will drop the superscripts writing $\be_i$ for
    $\be_i^\Pi$, etc. 
    \begin{enumerate}
    \item Obvious from the definition of $\be_j$. 
    \item To see that $\< \cdot, \cdot \>^\Pi$
        symmetrizes to $(\cdot, \cdot)$, we use the identity
    $\al_i=\be_i-\sum_{k\to i}\be_k$ (see Proposition~\ref{p:beta}),
    which    gives 
     \begin{equation*}
        \< \alpha_{i} ,\alpha_{j} \>^\Pi = \<
        \beta_{i} -\sum_{k\to i} \beta_{k} ,\alpha_{j}
        \>^\Pi = \begin{cases} 
                                 1 & \text{ if } i=j \\ 
                                -1 & \text{ if } j \to i \\ 
                                 0 & \text{otherwise}
                        \end{cases} 
    \end{equation*}
    Then
    \begin{align*}
    \< \alpha_{i} ,\alpha_{j} \>^\Pi + \< \alpha_{j},
        \alpha_{i} \>^\Pi &= 
                 \begin{cases} 2 & \text{if }i=j \\ 
                           -1 & \text{if }i\to j\: \text{or }j \to i \\ 
                            0 & \text{otherwise}
                \end{cases}\\
        &=(\alpha_{i} ,\alpha_{j}).
    \end{align*}

    To prove relation \eqref{e:serre_duality}, it suffices to prove
    that for all $k,i$ one has
    \begin{equation}\label{e:serre2}
    \<\be_k,C^{-1}\be_i\>^\Pi=-\<\be_i,\be_k\>^\Pi  .
    \end{equation}
    We will prove it by fixing $k$ and using  induction  in $i$, using the
    partial order defined  by the orientation. Thus, we can assume that
    \eqref{e:serre2} is true    for all $j\geqslant i$.

    Using  Proposition~\ref{p:Cbeta}, we can rewrite the
    left-hand side of \eqref{e:serre2} as 
    $$
    \<\be_k,C^{-1}\be_i\>^\Pi=-\<\be_k,\be_i\>^\Pi
                +\sum_{j \to i} \<\be_k,\be_j\>^\Pi
                +\sum_{j \leftarrow i}\<\be_k,
                                        C^{-1}\be_j\>^\Pi .
    $$
    The first two terms can be rewritten as
    $$-\<\be_k,\be_i\>^\Pi 
            +\sum_{j \to i}\<\be_k,\be_j\>^\Pi
        =-\<\be_k,\al_i\>^\Pi=-\delta_{ik} .
    $$
    The last term, using induction assumption, can be rewritten as 
    $$\sum_{j \leftarrow i} \<\be_k, C^{-1}\be_j\>^\Pi=
        -\sum_{j \leftarrow i}\<\be_j,\be_k\>^\Pi
       =-\sum_{j \leftarrow i}(\text{number of paths }j\to \cdots\to k).
    $$
    Thus, the left-hand side of \eqref{e:serre2} is 
    \begin{align*}
    \<\be_k,C^{-1}\be_i\>^\Pi&=-\delta_{ik}
     - \sum_{j \leftarrow i}(\text{number of paths }j\to \cdots \to k)\\
     &=-(\text{number of paths }i\to \cdots \to k)
        =-\<\be_i,\be_k\>^\Pi
    \end{align*}
    which proves \eqref{e:serre2}.
   
    \end{enumerate}
\end{proof}

\begin{thm}\label{t:form2}
    For fixed $C$, the form $\< \cdot,\cdot\>^\Pi$ does not
    depend on the choice of simple root system $\Pi$ compatible with
    $C$. Thus, we will denote this form       
    $\<\cdot,\cdot\>$ and call it the {\em Euler form}    defined
    by $C$.
\end{thm}
\begin{proof}
    Consider the difference $\ll \cdot , \cdot
    \gg=\<\cdot , \cdot\>^{\Pi_1}-\<\cdot ,
    \cdot\>^{\Pi_2}$. Since these two forms have the same
    symmetrization, $\ll \cdot , \cdot
    \gg$ is skew-symmetric and satisfies \eqref{e:serre2}. Thus, 
    $$
            \ll x,y \gg = -\ll y,C^{-1}x \gg =\ll C^{-1}x ,y \gg
     $$
     so $\ll (1-C^{-1} )x, y \gg =0$. Since $1$ is not an
        eigenvalue for $C^{-1}$ (see \cite[Lemma 8.1]{kostant1}), the operator $1-C^{-1}$ is invertible, so
        the form $\ll  \cdot ,\cdot \gg$ must be
        identically zero.
\end{proof}
\begin{prop} 
    Let $\< \cdot,\cdot\>$ be the Euler form defined
    in Theorem~\ref{t:form2}. Then 
    \begin{enumerate}
        \item The form $\< \cdot,\cdot\>$ is $C$--invariant:
            $\<Cx,Cy\>=\<x,y\>$. 
        \item
            $\< x,y\> =(x,(1-C^{-1})^{-1}y)=((1-C)^{-1}x,y)$
        
        \item Let $\omega_i^\Pi$ be fundamental weights:
        $(\omega_i^\Pi,\al_j^\Pi)=\delta_{ij}$. Then 
        $\beta_{i}^\Pi =(1-C) \omega^\Pi_{i}$
        
        \item The operator $(1-C)$ is an isomorphism $P(R)\to Q(R)$,
            where $Q(R)$, $P(R)$ are root and weight lattices of $R$
            respectively. 
    \end{enumerate}
\end{prop}

\begin{proof}\par\noindent
    \begin{enumerate}
    \item Choose a compatible simple root system $\Pi$; then by
        definition $\<\be_i^\Pi,\al_j^\Pi\>=\delta_{ij}$. On the other
        hand, $C(\Pi)$ is also a compatible simple root system, and
        $\al_i^{C(\Pi)}=C\al_i^\Pi$, $\be_i^{C(\Pi)}=C\be_i^\Pi$, so we
        get 
        $$
          \<C\be_i^\Pi,C\al_j^\Pi\>=\delta_{ij}=\<\be_i^\Pi,\al_j^\Pi\>.
        $$

    \item
        First we write $\< x,y \> =(x,Ay)$ for some $A$.
        Then using \eqref{e:serre_duality} we get:
        \begin{align*}
        (x,y) &= \< x,y \> +\< y,x \> = \< x,y \> -\< x, C^{-1}y \> &\\
              &=(x,Ay)-(x,AC^{-1} y) &\\
            &=(x,A(1-C^{-1}) y) . &
        \end{align*}
        Hence $A=(1-C^{-1})^{-1}$. (Note that $1-C^{-1}$ is invertible
        since 1 is not an eigenvalue for $C^{-1}$.) The second identity is
        proved in a similar way.

    \item From Part (2) we have that $\delta_{ij} =((1-C)^{-1}
    \beta_{i}^\Pi ,\alpha^\Pi_{j} )$ and hence $\beta_{i}^\Pi
    =(1-C) \omega_{i}^\Pi$.
    
    \item Since $\omega_i$ form a basis of $P(R)$, and $\be_i$ form a
    basis of $Q(R)$ (Proposition~\ref{p:beta}), this follows from part
    (3). 
\end{enumerate}
\end{proof}

As an immediate corollary, we get the following result, describing
the Euler form in $\Ihat$.

\begin{prop}\label{p:euler_ihat}
    There exists a unique function 
    $\< \cdot ,\cdot\>_{\Ihat} \colon \Ihat\times \Ihat\to\Z$
    satisfying 
    \begin{enumerate}

    \item $\< (i,n),(j,n) \>=\delta_{ij}$, \\
          $\< (i,n),(j,n+1) \>$ =the number of paths 
                $(i,n)\to \cdots \to  (j,n+1)$\\
            $=$number of edges between $i,j$ in $I$.

    \item For any $q=(k,m)\in \Ihat$ we have 
        $$ 
        \< q, (i,n) \>_{\Ihat} 
          -\sum_{j-i} \< q, (j,n+1) \>_{\Ihat} 
        +\< q, (i,n+2) \>_{\Ihat} = 0.
        $$

\end{enumerate}

\end{prop}
\begin{proof}
    Uniqueness easily follows by induction: for fixed $q=(i,n)$, 
    condition (1) defines $\<q, (*,n)\>_{\Ihat}$, $\<q,
    (*,n+1)\>_{\Ihat}$. Then condition (2) can be used to define
    $\<q, (*,n+2)\>_{\Ihat}$. Continuing in
    this way, we see that these two conditions completely determine
    $\<\cdot, \cdot\>_{\Ihat}$. 

    To prove existence, note that the form 
    $\<q_1,q_2\>_{\Ihat} =\<\Ph^{-1}(q_1), \Ph^{-1}(q_2)\>$, where $\< \cdot ,\cdot \>$ 
    is the Euler form on $R$ defined in Theorem~\ref{t:form2}, 
    satisfies all required properties. 
\end{proof}

Note that the proof of uniqueness actually gives a very simple and 
effective algorithm for computing $\<\cdot,\cdot\>_{\Ihat}$. 

\begin{thm}
    For a simply-laced Dynkin diagram $I$, define the set $\Ihat$ and
    lattice $Q$ by    \eqref{e:Ihat}, \eqref{e:Q} respectively. Let
    $\<\cdot, \cdot\>_{\Ihat}$ be the Euler form
    defined in Proposition~\ref{p:euler_ihat}, and let 
    $(x,y)_\Ihat=\<x,y\>_\Ihat+\<y,x\>_\Ihat$.

    Then $(\cdot,\cdot)_\Ihat$ is a positive definite
    symmetric    form on $Q$, and $\Ihat\subset Q$ is a root system with
    Dynkin    diagram $I$. 
\end{thm}
\begin{proof}
    Let $R$ be a root system with Dynkin diagram $I$, and $C$ a
    Coxeter element in the corresponding Weyl group. Then the map
    $\Phi$ constructed in Theorem~\ref{t:Phi2} identifies $R\to \Ihat$,
    $Q(R)\to Q$ and the inner product $(\cdot,\cdot)$ in $Q$
    with $(\cdot,\cdot)_\Ihat$.
\end{proof}

\section{Positive roots and the longest element in the Weyl group}\label{s:Delta}
    
Let $\Pi$ be a simple root system compatible with $C$ and $R_+^\Pi$ the corresponding set of positive roots. Then
bijection $\Phi\colon R\to \Ihat$ constructed in
Theorem~\ref{t:Phi2} identifies  $R_+^\Pi$ with a certain subset in
$\Ihat$. This subset can be identified with the usual Auslander--Reiten
quiver of the category of representations of quiver $I$ with the orientation defined by $\Pi$; this will
be discussed in detail in Section~\ref{s:identification}. 

In this section, we give an explicit description of the set
$\Phi(R_+^\Pi)$ in terms of $\Ihat$. 

\begin{figure}[ht!]
\includegraphics{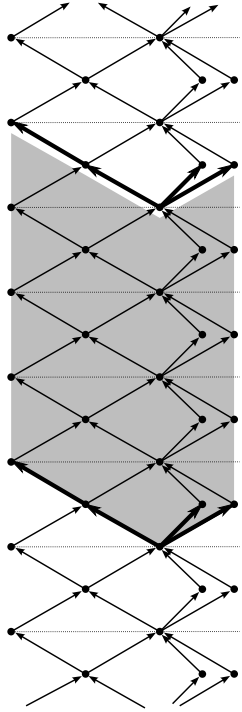}
\caption{Positive roots in $\Ihat$, for diagram of type $D_5$.  Bold
lines show $\Phi(\be_i^\Pi)$ and $\Phi(\be_i^{-\Pi})$; the shaded area
is the set $\Delta=\Phi(R_+^\Pi)$.}\label{f:Delta}
\end{figure}

\begin{thm}\label{t:w_0}
    Let $\Pi$  be a simple root system compatible with $C$, and
    $-\Pi=\{-\al\st \al\in \Pi\}$ the opposite simple root system. Let
    $h^\Pi$, $h^{-\Pi}\colon I\to \Z_{2h}$ be the corresponding
    height functions  as  defined in Theorem~\ref{t:height_functions}. 

    Let $\Delta^{\Pi}$ be the set of all vertices of $\Ihat$ ``between''
    correponding slices: 
    \begin{equation}
        \Delta^\Pi = \{(i,n)\in \Ihat\st h^\Pi(i)\le n< h^{-\Pi}(i)\}.
    \end{equation}
    \textup{(}see Figure~\ref{f:Delta}\textup{)}. We will consider 
    $\Delta^\Pi$ as a  quiver, with same edges as $\Ihat$. 

    \begin{enumerate}
         \item The map $\Phi$ gives identification $R_+^\Pi\to
            \Delta^\Pi$. 
         \item Define a partial order on $\Delta^{\Pi}$ by
            $q_1\preceq q_2$  if there exists a path from $q_1$ to $q_2$
            in $\Delta^{\Pi}$, and extend it to a complete order, writing  
            $$
                \Delta^{\Pi}=\{\al(1)=(i_1,n_1),\al(2)=(i_2,n_2),\dots,
                \al(l)=(i_l,n_l)\} 
            $$
            so that  $\al(a)\preceq\al(b)\implies a\leq b$. Then 
            $$
            s^\Pi_{i_1}\dots s_{i_l}^\Pi
            $$
            is a reduced expression for the longest element $w_0^\Pi$ of
            the  Weyl group. 
    \end{enumerate}
\end{thm}

\begin{proof}
    First, we prove that $\Delta^\Pi\subset \Phi(R_+^\Pi)$. Indeed,
choose    $i\in I$ and let $k$ to be the smallest positive integer such
that    $\be =\Phi^{-1}(i,h^\Pi(i)+2k)\notin R_+^\Pi$. Then $\be$
satisfies    $\be\in R_+^{-\Pi}$, $C^{-1}\be\in R_+^\Pi$. By
    Proposition~\ref{p:beta} applied to $-\Pi$, we see that we must
    have $\be=\be_i^{-\Pi}$, so $h^\Pi(i)+2k=h^{-\Pi}(i)$. Therefore,
    for all $n$ satisfying $h^\Pi(i)\le n< h^{-\Pi}(i)$, we have
    $(i,n)\in \Phi(R_+^\Pi)$. 

    Reversing the roles of $\Pi,-\Pi$, we see that the set
    $\Delta^{-\Pi}=\{(i,n)\in \Ihat\st h^{-\Pi}(i)\le n< h^{\Pi}(i)\}$
    staisfies $\Delta^{-\Pi}\subset\Phi(R_+^{-\Pi})=\Phi(-R_+^\Pi)$.
    Since $\Ihat=\Delta^{\Pi}\sqcup\Delta^{-\Pi}$, we must have
    $\Delta^\Pi=\Phi(R_+^\Pi)$. In particular,
    $l=|\Delta^\Pi|=|R_+^\Pi|$.

    To prove the second part, consider a sequence of sets
    $\Delta_0=\Delta^\Pi$,
    $\Delta_1=\Delta_0\setminus\{\al(1)\}=\{\al(2),\dots,\al(l)\}$, 
 $\Delta_k=\Delta_{k-1}\setminus\{\al(k)\}=\{\al(k+1),\dots,\al(l)\}$, 
 $\Delta_l=\varnothing\}$. 

    It is immediate from the definitions that one can
    write 
    \begin{equation}
            \Delta_k = \{(i,n)\in \Ihat\st h_k(i)\le n< h^{-\Pi}(i)\}.
    \end{equation}
    for the sequence of height functions $h_0=h^{\Pi}$,
    $h_1=s_{i_1}h_0$, \dots, $h_k=s_{i_k}h_{i_{k-1}}$, \dots,
    $h_l=h^{-\Pi}$, and  $i_k$ is a source for
    $h_{k-1}$.    
    In the same way, we can
    define a sequence    of corresponding simple root    systems       
    $\Pi_0=\Pi$,    $\Pi_1=s_{i_1}\Pi, \dots, \Pi_l=-\Pi$. One easily
    sees that
    $$
    \Pi_2=s_{i_2}^{\Pi_1}(\Pi_1)=s_{i_1}^\Pi s_{i_2}^\Pi s_{i_1}^\Pi
    (\Pi_1)=s_{i_1}^\Pi s_{i_2}^\Pi(\Pi).
    $$
    Repeating this argument, we get that $\Pi_k=s^\Pi_{i_1}\dots
    s_{i_k}^\Pi(\Pi)$; in particular, 
    $$
    -\Pi=\Pi_l=s_{i_1}^\Pi\dots s_{i_l}^\Pi(\Pi)
    $$
    Since $l=|R_+|=l(w_0)$, we see that $s_{i_1}^\Pi\dots s_{i_l}^\Pi$
    is a reduced expression for the longest element $w_0^\Pi$. 
\end{proof}

\begin{example}
Consider the case $R=A_{4}$ from Example~\ref{ex:phi-A4}. We obtain the expression $w_{0} = s_{1}s_{2}s_{1}s_{3}s_{2}s_{4}s_{1}s_{3}s_{2}s_{1}$ for the longest element.
\end{example}

Note that the second part of the theorem is equivalent to the algorithm
for constructing a reduced expression for $w_0$ in terms of the
Auslander--Reiten quiver, given in \cite{bedard} (in the form given here
it is reformulated in \cite[Theorem~1.1]{zelikson}). 

\begin{example}
Let $\Pi$ be such that it defines a bipartite orientation on $I$, as in 
Example~\ref{ex:bipartite}, so that 
$C=(\Pi_{i\in I_{0}} s_{i}) (\Pi_{i\in I_{1}} s_{i})$. Assume additionally that 
the Coxeter number is even: $h=2g$. Then it is known that $w_0=C^g$
(see \cite{kostant}), and thus 
$-\Pi=C^g(\Pi)$, and the corresponding height functions are related by 
$h^{-\Pi}=h^{\Pi}+2g$. In this case, the set $\Phi(R_+^\Pi)$ is a
``horizontal strip''.  
However, as the example of type $A$ shows (see Example~\ref{ex:phi-A4}), 
in general the set $\Phi(R_+^\Pi)$ can have a more complicated shape. 
\end{example}

For any $i\in I$ we define $i \ \check{}$ by $-\alpha_{i}^{\Pi} = w_{0}^{\Pi} (\alpha_{i \ \check{}}^{\Pi})$, where $w_{0}^{\Pi} \in W$ is the longest element. Equivalently we can define $i \ \check{}$ by $-\beta_{i}^{\Pi} = \beta_{i \ \check{}}^{-\Pi}$, where $\beta_{i}^{\Pi}$ are the $C$-orbit representatives defined above. Thus for the root systems of type $A, D_{2n+1}, E_{6}$ this map corresponds to the diagram automorphism, while for $D_{2n}, E_{7}, E_{8}$ this map is just the identity (corresponding to the fact that $-Id=C^{h/2} \in W$).

Now define a map $\nu_{\Ihat} :\Ihat \to \Ihat$ by the formula 
\begin{equation}\label{e:nuihat}
\nu_{\Ihat} (i,k) = (i \ \check{}, k+h)
\end{equation}
This map will correspond to the Nakayama permutation defined in Section~\ref{s:RepI}.

To see this is well defined we need to check that the image does in fact lie in $\Ihat$.  Note that if h is even $p(i \ \check{}) = p(i)$ and $k+h=k \mod 2$, so $(i \ \check{}, k+h) \in \Ihat$. If h is odd, then $R=A_{2n}$, $h=2n+1$ and $i \ \check{} = 2n-i+1$, so that $p(i \ \check{})=p(i) + 1$ and $k+h=k+1 \mod 2$, so  $(i \ \check{}, k+h) \in \Ihat$. Hence the map $\nu_{\Ihat}$ is well-defined.

\begin{lemma}\label{l:nu}
Under the identification $\Phi :R\to \Ihat$ the map $-Id$ corresponds to $\nu_{\Ihat}$.
\end{lemma}

\begin{proof} 
$\Phi (-\beta_{i}^{\Pi}) = \Phi (\beta_{i \ \check{}}^{-\Pi}) = (i \ \check{} , h^{-\Pi} (i \ \check{})) = \nu_{\Ihat} (i, h^{\Pi})$.

For the last equality, since the maps $\Phi$ and $\nu_{\Ihat}$ are compatible with the simple reflection $s_{j}$ when $j\in I$ is a sink or source, it is enough to consider the case where  $\Pi$ gives a bipartite splitting $I=I_{0} \sqcup I_{1}$ and $C=(\prod_{i\in I_{0}} s_{i}^{\Pi})(\prod_{i\in I_{1}} s_{i}^{\Pi})$. 

If $h=2g$ is even, then $w_{0}^{\Pi} = C^{g}$ (see \cite{kostant}) and hence $h^{\Pi}(i) + h = h^{-\Pi}(i \ \check{})$. 

If $h=2g+1$ is odd, then $R=A_{2g}$. We have that $w_{0}^{\Pi} = (\prod_{i\in I_{0}} s_{i}^{\Pi}) C^{g}$ and $i \ \check{} = 2g+1-i$, so that $i \ \check{} \in I_{1}$ if $i\in I_{0}$ and $i \ \check{} \in I_{0}$ if $i\in I_{1}$.\
For $i \ \check{} \in I_{0}$ we have that 
$$h^{-\Pi}(i \ \check{}) = w_{0}^{\Pi} h^{\Pi} (i \ \check{}) = h^{\Pi}(i \ \check{}) +2g+2 = h^{\Pi} (i) -1 +2g+2 = h^{\Pi} (i) +h.$$
For $i \ \check{} \in I_{1}$ a similar calculation shows that $h^{-\Pi}(i \ \check{}) = h^{\Pi} (i) + h$.

\end{proof}

Recall that a choice of compatible simple roots $\Pi$ gives a height function $h^{\Pi}$, and hence a slice given by $I_{\Pi} =\{ (i, h^{\Pi}(i)) \}$. The following Proposition gives another characterisation of $\Delta^{\Pi} \subset \Ihat$.

\begin{prop}\label{p:Delta}
 $\nu_{\Ihat} (I_{\Pi}) = I_{-\Pi}$, and so the subset $\Delta^{\Pi} \subset \Ihat$ can be identified as the subquiver lying between $I_{\Pi}$ and $\nu_{\Ihat} (I_{\Pi}) = I_{-\Pi}$.
\end{prop}
\begin{proof}
This follows immediately from Lemma~\ref{l:nu} and the identification of $\Delta^{\Pi}$ with the set of roots lying between $h^{\Pi}$ and $h^{-\Pi}$ given in Theorem~\ref{t:w_0}.

\end{proof}

\section{The Auslander-Reiten quiver of $\Rep(I_{\Om})$ and $\Z I$}\label{s:RepI}

We now relate the previous sections to the theory of quiver representations.
In this section we recall the definition of the Auslander-Reiten quiver, as well as some basic facts about quiver representations.

\begin{figure}[ht]
	\includegraphics[height=3.00in]{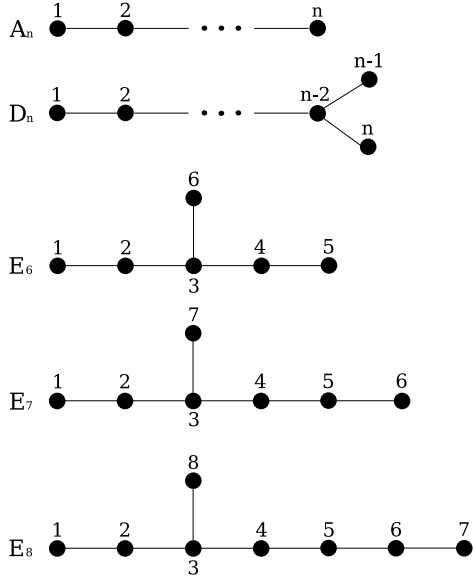}
	\caption{Identification of vertices of Dynkin Diagram with $\{ 1, \ldots, r \}$.}\label{f:dynkindiagrams}
\end{figure}

Let $I$ be a Dynkin diagram of type $A,D,E$. Choose an identification of $I$ with $\{ 1, \ldots , r \}$ as in Figure~\ref{f:dynkindiagrams}.

Given an orientation $\Om$ on $I$, consider the category $\Rep(I,\Om)$ of representations of the quiver $(I,\Om)$. For a representation $X$, let $X_{i}$ denote the vector space at the vertex $i\in I$. 
The dimension vector of a representation $X$ is given by $\dim X = (\dim X_{1}, \ldots, \dim X_{r} ) \in \Z^{I} $. This gives a map $\dim : K(\Rep(I, \Om)) \to \Z^{I}$.

Let $\Ind(I,\Om)$ be the set of isomorphism classes of indecomposable objects in $\Rep(I,\Om)$. The Auslander-Reiten quiver $\Gamma (I, \Om)$ of $\Rep (I, \Om)$ has $\Ind(I, \Om)$ as its set of vertices, and has oriented edges given by $\left[ X \right] \to \left[ Y \right]$ if there exists an indecomposable morphism $f: X \to Y$.

In what follows the projective representations in $\Ind (I, \Om)$ play an important role, so we will review the description of the projectives given in \cite{gabriel}. For $i \in I$, we define the representation $P(i)$ by 
\begin{equation}\label{e:proj}
P_{j}(i) = \Path_{(I, \Om)} (i,j) = \begin{cases} \CC & \text{ if there is a path i}  \to \cdots \to \text{j in} \ (I, \Om) \\ 0 & \text{otherwise} \end{cases}
\end{equation}
Then the representations of the form $P(i)$ form a complete list of representatives of projectives in $\Ind(I, \Om)$ (see \cite{gabriel} Section 4.4 for the proof).

Given $I$ we construct a quiver $\Z I$ following \cite{gabriel}. Using the identification of $I = \{ 1, \ldots , r \}$ we construct $\Z I$ as follows. The vertices are pairs $(i,k)$ with $i\in I$ and $k\in \Z$. For $i,j$ connected in $I$ and $i < j$, we join $(i,k)$ to $(j,k)$ and $(j,k)$ to $(j,k+1)$.

On $\Z I$ the Auslander-Reiten translation $\tau :\Z I \to \Z I$ is given by 
\begin{equation}\label{e:ARtrans}
\tau (i,k) = (i,k-1)
\end{equation}

\begin{example}
The quiver $\Z A_{4}$ is depicted in Figure~\ref{f:ZA4} 
\begin{figure}[ht]
	\includegraphics[height=3.00in]{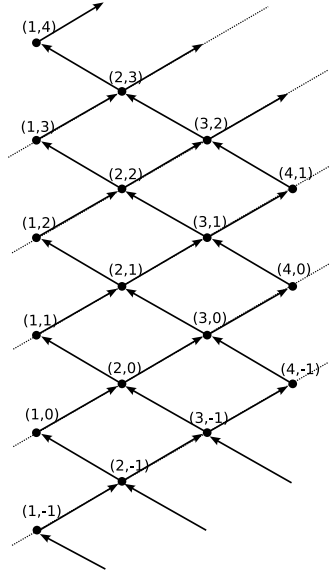}
	\caption{The quiver $\Z A_{4}$.}
\label{f:ZA4}
\end{figure}
\end{example}

\begin{remark}
Note that here the $\Z$ direction is vertical, while in \cite{gabriel} and most of the literature the $\Z$ direction is horizontal. This choice was made to agree with the $\Z_{2h}$ direction in $\Ihat$. 
\end{remark}

On $\Z I$  the Nakayama permutation $\nu :\Z I \to \Z I$ is defined by the formulas below (see \cite{gabriel}):
\begin{equation}\label{e:nu}
\nu (i,k) = (i \ \check{} ,k+i ) \ \text{for type} \ A
\end{equation}
\begin{equation}
\nu(i,k)= (i \ \check{} , k+h/2) \ \text{for types} \  D, E
\end{equation}

\begin{remark}
Note that we have shifted the Nakayama permutation $\nu$ given in \cite{gabriel} by $\tau^{-1}$, so that the map will agree with $\nu_{\Ihat}$.
\end{remark}

Following \cite{gabriel}, a slice of $\Z I$ is a connected full subquiver of $\Z I$ which contains a unique representative of $\{ (i,n) | n\in \Z \}$ for each $i\in I$.

Identify $1\in I$ with $(1,0) \in \Z I$. For each orientation $\Om$ of $I$, this identifies the quiver $(I,\Om)$ with a unique slice passing through $(1,0)$. We denote this slice by $I_{\Om}$.

The following Theorem describes the Auslander-Reiten quiver explicitly as a subquiver of $\Z I$. 
 
\begin{thm}\cite{gabriel}\label{t:AR}
The Auslander-Reiten quiver $\Gamma (I, \Om)$ of $\Rep(I, \Om)$ can be identified with the full subquiver of $\Z I$ lying between the slice $I_{\Om ^{opp}}$ and the slice $\nu (I_{\Om^{opp}} )$. Morever, the projective representations $P_{i}$, correspond to the slice $I_{\Om^{opp}}$.
\end{thm}
For a proof see \cite[Proposition, p.50]{gabriel}.

\begin{example}\label{ex:ARA4}
For the case $I =A_{4}$ and the orientation $\Om$ given by $1\leftarrow 2 \leftarrow 3 \leftarrow 4$ the Auslander-Reiten quiver is shown in Figure~\ref{f:ARA4} as the shaded region.
\begin{figure}[ht]
	\includegraphics[height=3.00in]{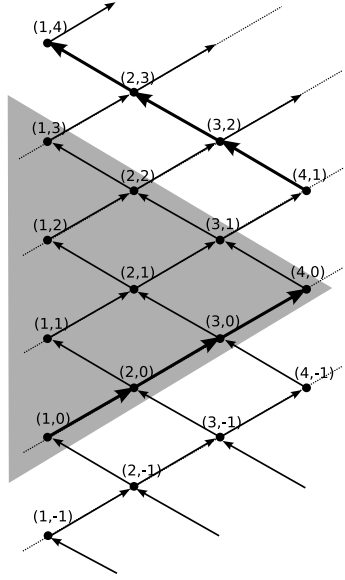}
	\caption{The Aulander-Reiten quiver is shown as the shaded region. The slices $I_{\Om^{opp}}$ and $\nu (I_{\Om^{opp}})$ are shown in bold.}
\label{f:ARA4}
\end{figure}
\end{example}

\section{Relationship between $\Z I$ and  $\Ihat$}\label{s:identification}

In this section we show that for any choice of $\Pi$ compatible with $C$ there is an identification of the Auslander-Reiten quiver with the subset $\Delta_{+}^{\Pi} \subset \Ihat$.

Let $\Pi$ be a simple root system compatible with $C$. This gives a height function $h^{\Pi}$, and hence gives an orientation $\Om$ on I, as well as a slice $I_{\Om}$. Define a covering map $P_{\Pi} : \Z I \to \Ihat$ as  follows:
\begin{equation}\label{e:covermap}
P_{\Pi}(i,k) = (i, h^{\Pi} (1) + l(i) +2k)
\end{equation}
where $l(i)$ is the number of edges in $I$ between $1$ and $i$. 

To see that this preserves arrows, note that if $(i,k) \to (j,n)$ in $\Z I$, then $i$ is connected to $j$ in $I$ and either $k=n$ or $k=n-1$ and we see that $P_{\Pi} (i,k) \to P_{\Pi} (j,n) $ in $\Ihat$.

\begin{prop}\label{p:P_Pi}
\par\indent
\begin{enumerate}
\item $P_{\Pi} \circ \tau = \tau_{\Ihat}^{-1} \circ P_{\Pi}$.
\item $P_{\Pi} \circ \nu = \nu_{\Ihat} \circ P_{\Pi}$
\item The map $P_{\Pi}$ identifies the Auslander-Reiten quiver $\Gamma (I, \Om^{opp})$ in $\Z I$ with $\Delta^{\Pi}$ in $\Ihat$.
\item Under the identification $P_{\Pi}$, the projective representation $P_{i}$ corresponds to $\beta_{i}^{\Pi}$.
\end{enumerate}
\end{prop}

\begin{proof}
\par\indent
\begin{enumerate}

\item \begin{align*}
P_{\Pi} (i, k-1) &= (i, h(1) + l(i) + 2(k-1) ) &\\  
&=(i, h(1) + l(i) + 2k -2) &\\ 
&= \tau_{\Ihat}^{-1} P_{\Pi} (i,k) &
\end{align*}
\item For $D,E$ type note that we have $\nu (i,k) = (i \ \check{} ,k+h/2)$ so that 
\begin{align*}
P_{\Pi}(\nu (i,k)) &= P_{\Pi} (i \ \check{}, k+h/2) &\\
&= (i \ \check{}, h(1) + l(i \ \check{}) + 2(k+h/2) ) &\\
&= (i \ \check{}, h(1) +l(i) +2k +h) &\\
&= \nu_{\Ihat} (i, h(1)+l(i) +2k) &\\
&= \nu_{\Ihat} \circ P_{\Pi} (i, k) &
\end{align*}

Similarly for $A$ type, we have that $h=n+1$, $l(i) = i-1$ and $\check{\imath} = n-i+1$. Then  we have

\begin{align*}
P_{\Pi} (\nu (i,k)) &= P_{\Pi} (i \ \check{} , k+i) &\\
&= (i \ \check{} , h(1) + l(i \ \check{}) +2(k+i) ) &\\
&= (i \ \check{}, h(1) + i \ \check{} -1 +2k +2i ) &\\
&= (i \ \check{}, h(1) +n-i+1 - 1 +2k+2i) &\\
&= (i \ \check{}, h(1) +i-1 +2k +n+1) &\\
&= (i \ \check{}, h(1) + l(i) +2k + h) &\\
&=\nu_{\Ihat} \circ P_{\Pi} (i,k) &
\end{align*}

\item Given $\Pi$ we obtain a unique height function $h^{\Pi}$, which also gives an orientation $\Om$, and a unique slice $I_{\Om}$ in $\Ihat$. If we consider $\Rep(I, \Om^{opp})$ we obtain a unique slice $I_{\Om}$ through $(1,0)$ in $\Z I$. By the construction of the map $P_{\Pi}$ these two slices are identified. By Part 2 we know that $P_{\Pi}$ identifies $\nu$ with $\nu_{\Ihat}$. Then the description of $\Gamma(I, \Om^{opp})$ given in Theorem~\ref{t:AR}, and the description of $\Delta^{\Pi}$ given in Proposition~\ref{p:Delta}, shows that the map $P_{\Pi}$ identifies  $\Gamma(I, \Om^{opp})$ with $\Delta^{\Pi}$. (Compare Example~\ref{ex:ARA4} with Example~\ref{ex:phi-A4}.)

\item The $\beta_{i}^{\Pi}$ map to the slice $I_{\Om} \subset \Ihat$. Similarly, the projective representations $P(i)$ map to the slice $I_{\Om} \subset \Z I$. By Part 3, these are identified by $P_{\Pi}$.
\end{enumerate}
\end{proof}

Given a simple root system $\Pi$, compatible with $C$, we have obtained a bijection of the Auslander-Reiten quiver $\Gamma (I, \Om^{opp})$ of the category $\Rep (I, \Om^{opp})$ with the subquiver $\Delta^{\Pi} \subset \Ihat$. Both $\Gamma (I,\Om^{opp})$  and $\Delta^{\Pi}$ correspond to the set of positive roots $R_{+}^{\Pi}$. For $\Gamma (I, \Om^{opp})$ the correspondence is the usual identification between $\Ind(I, \Om^{opp})$ and positive roots $R_{+}^{\Pi} \subset \oplus _{i\in I} \Z\alpha_{i}$, given by the dimension vector $\dim X $. For $\Delta^{\Pi}$ it is given by the map $\Phi :R\to \Ihat$ from Section~\ref{s:Ihat}. The following Theorem shows that the bijection $\Gamma (I, \Om^{opp}) \to \Delta^{\Pi}$ agrees with these identifications to $R_{+}^{\Pi}$.

\begin{thm}
Let $\Pi$ be a compatible simple root system, $h^{\Pi}$ the corresponding height function and $\Om$ the corresponding orientation. 
The following diagram is commutative: 
$$\xymatrix{
& \Gamma (I, \Om^{opp}) \ar[dl]_{\dim} \ar[dd]^{P_{\Pi}} \\
R_{+}^{\Pi} \ar[dr]_{\Phi} & \\
& \Delta^{\Pi} \\
}$$

\end{thm}
\begin{proof}
Under the map $\dim :\Gamma (I, \Om^{opp}) \to R_{+}^{\Pi} $ the projective representation $P(i)$ given by 
$$P_{j}(i) = \begin{cases} \CC & \text{j} \leqslant \ \text{i} \\ 0 & \text{otherwise} \end{cases}$$ is identified with the $C$-orbit representative $\beta_{i}^{\Pi}$. By Theorem~\ref{t:AR} these representations form the slice $I_{\Om} \subset Z I$. The $C$-orbit representatives map to the slice $I_{\Om}$ in $\Ihat$ under the map $\Phi$. By Proposition~\ref{p:P_Pi}  these two slices are identified by the map $P_{\Pi}$. Hence the diagram commutes on the projectives $P(i)$. The map $\dim$ identifies $\tau^{-1}$ (where $\tau$ is the Auslander-Reiten translation) with $C$, the map $\Phi$ identifies $C$ with $\tau_{\Ihat}$, and $P_{\Pi}$ identifies $\tau^{-1}$ with $\tau_{\Ihat}$. Since any $X\in \Gamma(I, \Om^{opp})$ is of the form $\tau^{-k} (P(i))$, and any $\alpha \in R_{+}^{\Pi}$ is of the form $C^{k} \beta_{i}^{\Pi}$ we see that the diagram commutes.
\end{proof}

\section{Identification of Euler Forms}\label{s:eulerform}
In this section we show that for a choice of compatible simple system $\Pi$, the bilinear form $\< \cdot , \cdot \>_{R}$ on $R$, constructed in Section~\ref{s:euler}, corresponds to the Euler form $\< X, Y \> = \dim \RHom (X,Y) = \dim \Hom (X,Y) - \dim \Ext^{1} (X,Y)$ on $\Rep (I, \Om^{opp})$, where $\Om$ is the orientation determined by $\Pi$.

Note that the Euler form satisfies
\begin{equation}\label{e:eulerform} 
\< X,Y \> = -\<Y, \tau X\> = \< \tau X , \tau Y\>
\end{equation}
(see \cite{crawley-boevey}, for a proof).

\begin{thm}
The map $P_{\Pi} :\Z I \to \Ihat$ identifies the Euler form $\< \cdot, \cdot \>$ on $\Gamma (I,\Om^{opp})$ with the form $\< \cdot , \cdot \>_{R}$ on $R_{+}^{\Pi}$.
\end{thm}
\begin{proof}
Recall that the $C$-orbit representative $\beta_{i}^{\Pi}$ corresponds to the projective representation $P(i)$, while the simple root $\alpha_{i}^{\Pi}$ corresponds to the simple representation $X(i)$. It is well known (see \cite{crawley-boevey} p.24) that $\< P(i), X(j) \> = \delta_{ij} $. \\
Define a form $\ll \cdot , \cdot \gg$ on $\Gamma (I, \Om)$ by $\ll X,Y \gg = \< P_{\Pi} (X) , P_{\Pi} (Y) \>_{R}$ where $\< \cdot , \cdot \>_{R}$ is the Euler Form on $R$ defined in Section~\ref{s:euler}. Then $\ll P(i),X(j) \gg = \< \beta_{i}^{\Pi} , \alpha_{j}^{\Pi} \>_{R} = \delta _{ij}$, and $\ll \cdot ,\cdot \gg$ satisfies Equation~\ref{e:eulerform}.
However,  since the value of $\ll P(i),X(j) \gg$ and Equation~\ref{e:eulerform} completely determine the form, the two forms $\< \cdot ,\cdot \>$ and $\ll \cdot , \cdot \gg$ are equal.
\end{proof}

\bibliographystyle{amsalpha}

\end{document}